\documentclass[11pt]{amsart}
\usepackage{amsmath,amssymb,amsthm}
\usepackage[latin1]{inputenc}
\usepackage{version,tabularx,multicol}
\usepackage{graphicx,float}

\newcommand{\reff}[1]{(\ref{#1})}

\theoremstyle{plain}
\newtheorem{theo}{Theorem}[section]

\newtheorem{cor}[theo]{Corollary}
\newtheorem{prop}[theo]{Proposition}

\newtheorem{lem}[theo]{Lemma}

\theoremstyle{remark}
\newtheorem{rem}[theo]{Remark}

\newcommand{\ca}{{\mathcal A}}
\newcommand{\cb}{{\mathcal B}}
\newcommand{\cc}{{\mathcal C}}

\newcommand{\cf}{{\mathcal F}}

\newcommand{\ci}{{\mathcal I}}

\newcommand{\cm}{{\mathcal M}}

\newcommand{\E}{{\mathbb E}}

\newcommand{\N}{{\mathbb N}}
\renewcommand{\P}{{\mathbb P}}

\newcommand{\R}{{\mathbb R}}

\newcommand{\ind}{{\bf 1}}

\newcommand{\Supp}{{\rm Supp}\;}

\newcommand{\val}[1]{\mathop{\left| #1 \right|}\nolimits}
\newcommand{\inv}[1]{\mathop{\frac{1}{ #1}}\nolimits}
\newcommand{\expp}[1]{\mathop {\mathrm{e}^{ #1}}}

\begin{document}

\title[Properties of the exploration process]{Feller  property and
  infinitesimal  generator of  the exploration process}

\date{\today}
  
\author{Romain Abraham} 

\address{
Romain Abraham,
MAPMO, Université d'Orléans,
B.P. 6759,
45067 Orléans cedex 2
France
}
  
\email{romain.abraham@univ-orleans.fr} 

\author{Jean-François Delmas}

\address{
Jean-Fran\c cois Delmas,
ENPC-CERMICS, 6-8 av. Blaise Pascal,
  Champs-sur-Marne, 77455 Marne La Vallée, France.}

\email{delmas@cermics.enpc.fr}

\thanks{The research of the second author was partially supported by 
  NSERC Discovery Grants of the Probability group at Univ. of British Columbia}

\begin{abstract}
  We  consider the exploration process associated to the continuous random tree (CRT)
  built using  a Lévy  process with  no negative
  jumps. This process has been studied by Duquesne, Le~Gall and
  Le~Jan. This measure-valued Markov process is a useful tool to study
  CRT as well as  super-Brownian motion with general branching
  mechanism. In this paper we prove this process is Feller, and we
  compute its infinitesimal generator on exponential functionals and
  give the corresponding martingale. 
  \end{abstract}



\maketitle


\section{Introduction}

The coding  of a  Lévy continuous random  tree (CRT) by  its exploration
process  can be found  in Le  Gall and  Le Jan  \cite{lglj:bplpep}. They
associate to a  Lévy process with no negative jumps  that does not drift
to $+\infty $, $X=(X_t, t\geq  0)$, a critical or sub-critical continuous
state branching process  (CSBP) and a Lévy CRT which  keeps track of the
genealogy of the CSBP.  The exploration process $(\rho_t,t\geq 0)$ takes
values  in the set  of finite  measures on  $\R_+$.  Informally,  for an
individual $t\geq 0$, its generation is recorded by its height $H_t=\sup
\Supp \rho_t$,  where $\Supp  \mu$ stand for  the closed support  of the
measure  $\mu$ on  $\R$.  And $\rho_t(dv)$  records,  for the  individual
labeled $t$,  the ``number''  of brothers of  its ancestor at  height $v$
whose labels are larger than $t$.

For  instance,  the height  process  $(H_t,  t\geq  0)$ of  Aldous'  CRT
\cite{a:crt3} is  a normalized  Brownian excursion, and  the exploration
process $(\rho_t, t\geq 0)$ is the Lebesgue measure, that is $\rho_t$ is
the Lebesgue measure on $[0,H_t]$.

The height process, which is enough to code a CRT, is not Markov whereas
the  exploration process is  Markov.  The  exploration process  has been
first introduced  in order to  study super-Brownian motion  with general
branching       mechanism        (see       \cite{dlg:rtlpsbp}       and
\cite{lglj:bplplfss}).  The strong  Markov property  of  the exploration
process   is   a  fundamental   property   and   has   been  proved   in
\cite{lglj:bplpep}. In  \cite{ad:falp}, the exploration  process is used
in order  to construct  a fragmentation process  associated to  a L\'evy
process  using Markov properties and martingale problems.

The  goal of  this paper  is  to give  some martingales  related to  the
exploration  process.  We  also   prove  that  the  exploration  process
satisfies  the Feller  property  (and hence  has  c\`ad-l\`ag paths  and
satisfies the strong  Markov property). And we compute  the generator of
the exploration process for  exponential functionals. At this point, let
us  mention  that  we  will  always  work  with  the  topology  of  weak
convergence for finite  measures. But, as the set  of finite measures on
$\R_+$ is  not locally  compact for this  topology, the  Feller property
does  not  make sense  on  that  space. So,  we  must  first extend  the
definition of the transition function  of the exploration process to the
set of finite measures  on $E=\R_+\cup\{+\infty\}$ endowed with a metric
that makes it compact.

We   give  the   general  organization   of  the   paper.    In  Section
\ref{sec:main}  we  recall  the   construction  and  definition  of  the
exploration process  and some of its  properties (which can  be found in
\cite{dlg:rtlpsbp}) that will be used  in the sequel.  We then state the
Feller   property,   Theorem   \ref{th:feller}.    We  also   give   its
infinitesimal   generator   for   exponential   functionals,   Corollary
\ref{cor:inf-gene}, as  well as  the related martingales  in Corollaries
\ref{cor:mart1}  and  \ref{cor:mart-ts}.    We  also  recall  in  Remark
\ref{rem:AD}   the  relation   given  in   \cite{ad:falp}   between  the
infinitesimal  generator of  the exploration  process associated  to the
Lévy  process  with  Laplace   exponent  $\psi$  and  the  infinitesimal
generator of the exploration process associated to the Lévy process with
Laplace    exponent   $\psi^{(\theta)}=\psi(\theta+\cdot)   -\psi(\theta)$,
$\theta\geq 0$.  Section \ref{sec:Feller} is devoted to the proof of the
Feller property.   We eventually compute the  infinitesimal generator of
the   exploration  process  for   exponential  functionals   in  Section
\ref{sec:inf-gene}.

\section{Definitions and main results}
\label{sec:main}

\subsection{Notations and exploration process}\label{sec:levysnake}

Let $S$  be a  metric space, and  $\cb(S)$ its Borel  $\sigma$-field. We
denote  by   $\cb_+(S)$  (resp.  $\cb_b(S)$)  the   set  of  real-valued
non-negative (resp.  bounded) measurable  functions defined on  $S$. Let
$\cm_f(S)$ be the set of finite measures defined on $S$ endowed with the
topology of weak convergence.   We write $\cm_f$ for $\cm_f(\R_+)$.  For
$\mu\in \cm_f(S)$, $f\in \cb_+(S)$  we write $\langle \mu,f \rangle$ for
$\int_{S} f(x) \; \mu(dx)$.

We consider a $\R$-valued Lévy process $X=(X_t,t\ge 0)$ with no negative
jumps, no Brownian part and  starting from 0.  Its law is characterized
by its Laplace transform: for $\lambda\ge 0$
\[
\E\left[\expp{-\lambda X_t}\right]=\expp{t\psi(\lambda)},
\]
and we suppose that its Laplace exponent, $\psi$, is given by 
\[
\psi(\lambda)=\alpha_0\lambda+\int_{(0,+\infty)}\pi(d\ell)
\left[\expp{-\lambda\ell}-1+\lambda\ell\right],  
\]
with  $\alpha_0\ge  0$  and  the   Lévy  measure  $\pi$  is  a  positive
$\sigma$-finite  measure  on   $(0,+\infty)$  such  that  $\displaystyle
\int_{(0,+\infty)}     (\ell\wedge     \ell^2)\pi(d\ell)<\infty$     and
$\displaystyle \int_{(0,1)}\ell\pi(d\ell)=\infty$.  The first assumption
(with the  condition $\alpha_0\ge 0$)  implies the process $X$  does not
drift to  $+\infty$, while the second  implies $X$ is  a.s.  of infinite
variation.   

For $\mu\in \cm_f$, we define 
\begin{equation}
   \label{eq:def-H-mu}
H^\mu=\sup \Supp (\mu)\in [0,\infty ],
\end{equation}
where $ \Supp(\mu)$ is the closed support of $\mu$,  with the
convention $\sup\emptyset=0$.

We recall the definition and properties of the exploration process which
are    given    in    \cite{lglj:bplpep},    \cite{lglj:bplplfss}    and
\cite{dlg:rtlpsbp}. The results of  this section are mainly extract from
\cite{dlg:rtlpsbp}.
 
Let
$I=(I_t,t\ge  0)$ be the  infimum process  of $X$:  $I_t=\inf_{0\le s\le
  t}X_s$.   We will also consider for every $0\le s\le
t$ the infimum of $X$ over $[s,t]$:
\[
I_t^s=\inf_{s\le r\le t}X_r.
\]
There exists a sequence $(\varepsilon_n,n\in \N^*)$ of positive real
numbers decreasing to 0 s.t. 
\[
\tilde H_t= \lim_{k\rightarrow\infty } \inv{\varepsilon_k} \int_0^t
\ind_{\{X_s<I^s_t+\varepsilon_k\}}\; ds
\]
exists and is finite a.s. for all $t\geq 0$. 

The point  0 is regular  for the Markov  process $X-I$, and $-I$  is the
local time of  $X-I$ at 0 (see \cite{b:pl}, chap. VII).  Let $\N$ be the
associated  excursion  measure  of  the  process $X-I$  out  of  0,  and
$\sigma=\inf\{t>0; X_t-I_t=0\}$ be the  length of the excursion of $X-I$
under $\N$. Recall that under $\N$, $X_0=I_0=0$.

{F}rom Section  1.2 in \cite{dlg:rtlpsbp},  there exists a  $\cm_f$-valued
process, $\rho^0=(\rho^0_t,  t\geq 0)$, called  the exploration process,
such that :
\begin{itemize}
\item   For  each   $t\geq   0$,  a.s.   $H^0_t=   \tilde  H_t$,   where
  $H_s^0=H^{\rho^0_s}$.
   \item For every $f\in \cb_+(\R_+) $,
\[
\langle \rho^0_t,f\rangle =\int_{[0,t]} d_sI_t^sf(H^0_s),
\]
or equivalently 
\begin{equation}\label{eq:def_rho}
\rho^0_t(dr)=\sum_{\stackrel{0<s\le t}
  {X_{s-}<I_t^s}}(I_t^s-X_{s-})\delta_{H_s^0}(dr). 
\end{equation}
\item Almost surely, for every $t>0$, we have 
$\langle \rho_t^0,1\rangle =X_t-I_t$. 
\end{itemize}

\subsection{The Feller property}

In  Proposition 1.2.3  \cite{dlg:rtlpsbp}, Duquesne  and  Le~Gall proved
that the  exploration process  is a strong  Markov process  with càd-làg
paths in $\cm_f$ equipped with  the topology of weak convergence (in fact
they prove  the exploration process  is a.s. càd-làg with  the variation
distance on  finite measures). Here,  we improve this result  by proving
that  the exploration  process fulfills  the Feller  property.   

We  set $E=\R_+\cup\{+\infty\}$ equipped  with a  distance that  makes $E$
compact. The set $\cm_f(E)$ is locally compact.

In the definition  of the exploration process, as $X$  starts from 0, we
have  $\rho_0=0$ a.s.  To  get the  Markov property  of $\rho$,  we must
define  the  process  $\rho$  started  at any  initial  measure  $\mu\in
\cm_f(E)$.  

Let $\mu\in \cm_f(E)$. For  $a\in [0,  \langle \mu,1\rangle ]  $, we  define the
erased measure $k_a\mu$ by
\[
k_a\mu([0,r])=\mu([0,r])\wedge (\langle \mu,1\rangle -a), \quad
\text{for $r\in E$}. 
\]
If $a> \langle  \mu,1\rangle $, we set $k_a\mu=0$.   In other words, the
measure $k_a\mu$ is the measure $\mu$ erased by a mass $a$ backward from
$H^\mu$, with the natural extension of \reff{eq:def-H-mu} to $\cm_f(E)$.

For $\nu,\mu \in \cm_f(E)$,  we define
the concatenation $[\mu,\nu]\in \cm_f(E) $ of the two measures, using the
convention $x+\infty =+\infty $ for $x\in E$, by
\[
\bigl\langle [\mu,\nu],f\bigr\rangle =\bigl\langle \mu,f\bigr\rangle
+\bigl\langle \nu,f(H^\mu+\cdot)\bigr\rangle , 
\quad f\in \cb_+(E).
\]

Eventually,  we  set  for  every  $\mu\in  \cm_f(E)$  and  every  $t>0$,
\[
\rho_t=\bigl[k_{-I_t}\mu,\rho_t^0].
\]
We say that  $\rho=(\rho_t, t\geq
0)$  is  the exploration  process  started  at  $\rho_0=\mu$, and  write
$\P_\mu$ for its law.  We set $H_t=H^{\rho_t}$. The process $\rho$ is an
homogeneous Markov process (this  is a direct consequence of Proposition
1.2.3 in \cite{dlg:rtlpsbp}).

\begin{theo}
   \label{th:feller}
The exploration process, $\rho$,  defined on $\cm_f(E)$ equipped with the
topology of weak convergence is Feller.
\end{theo}

The proof of  this theorem is given in  Section \ref{sec:Feller}.  Since
$\rho$ has  the Feller property we  recover it enjoys  the strong Markov
property and  has a.s. càd-làg  paths. Notice that if  $\mu\in\cm_f$, as
$\rho^0$  is $\cm_f$  valued,  we  get by  construction  that $\rho$  is
$\cm_f$  valued. We  recover  that the  exploration  process in  $\cm_f$
enjoys the strong  Markov property and has a.s.   càd-làg paths.  We can
consider the infinitesimal generator of $\rho$.

\subsection{The infinitesimal generator}

Let  $f$ a  bounded  non-negative  function defined  on  $\R_+$ of  class
$\cc^1$     with    bounded  first   derivative    such     that    $f(\infty
)=\lim_{t\rightarrow\infty       }       f(t)$   exists.

We  consider  $F,  K\in
\cb(\cm_f(E))$ defined by: for $\mu\in \cm_f(E)$,
\begin{equation}
   \label{eq:def-FK}
   F(\mu)=\expp{- \langle \mu,f \rangle},\quad\text{and}\quad
   K(\mu)=F(\mu) \left[\psi( f(H^\mu)) - 
     f'(H^\mu) \ind_{\{ H^\mu<\infty \}} \right] . 
\end{equation}
Notice that  $F$  is a continuous function,
whereas $K$ belongs only to $\cb(\cm_f(E))$ a priori.

We denote  by $P_t$ the  transition function of the  exploration process
$\rho$  at time  $t$ and  we  consider its  resolvent: for  $\lambda>0$,
$U_\lambda=\int_0^\infty \expp{-\lambda t} P_t$.

\begin{theo}
   \label{th:inf-gene-0}
Let $\lambda>0$.
We have 
\[
U_\lambda (\lambda F - K)(\mu)=\expp{- \langle \mu,f \rangle} -
\frac{f(0)}{\psi^{-1}(\lambda)} \expp{- \psi^{-1}(\lambda) \langle \mu,1
  \rangle}. 
\]
\end{theo}

The proof of this Theorem is given in Section \ref{sec:inf-gene}.  There
exists a local time at $\mu=0$ which has to be considered when computing
the  infinitesimal  generator.  This  translates  into  the
condition $f(0)=0$ for the domain of the infinitesimal generator, as
stated in the next Corollary.

\begin{cor}
   \label{cor:inf-gene}
Assume that $f(0)=0$, $\lambda>0$.
We have $U_\lambda (\lambda F - K)=F$. In particular $(F,K)$ belongs to
the extended domain of the infinitesimal generator of $\rho$. 
\end{cor}

\begin{rem}
 Notice the function $K$ is not continuous. However, it can be proved
 directly that if $\mu\neq 0$ with compact support then
\[
\lim_{\lambda\rightarrow      \infty } \lambda\Big( \lambda U_\lambda F
-F\Big)(\mu)=K(\mu)
\]
and if $\mu=0$ 
\[
\lim_{\lambda\rightarrow
     \infty } \lambda\Big( \lambda U_\lambda F -F+
   \frac{f(0)}{\psi^{-1}(\lambda)} \Big) (0)=\psi(f(0))- f'(0).
\]
\end{rem}

Standard results on Markov process implies the following Corollary, see
e.g. \cite{ek:mp} Chapter 4.

\begin{cor}
\label{cor:mart1}
 Assume $f(0)=0$, $\lambda\geq 0$. The process $(M_t, t\geq 0)$ defined by 
\[
M_t=\expp{- \lambda t - \langle \rho_t,f \rangle} + \int_0^t \expp{-
  \lambda s - \langle
  \rho_s,f \rangle} \left[\lambda -\psi( f(H_s)) + 
     f'(H_s)\ind_{\{H_s<\infty \}} \right]\; ds,
\]
where $H_s=H^{\rho_s}$, is a martingale w.r.t. the filtration generated
by $\rho$.  
\end{cor}

We deduce from the optional stopping Theorem and the monotone
class Theorem the next result. (Notice we don't assume $f(0)=0$.)

\begin{cor}
\label{cor:mart-ts}
 The process $(M_t, t\geq 0)$ defined by 
\[
M_t=\expp{-  \langle \rho_{t\wedge \sigma} ,f \rangle} - \int_0^{t\wedge
  \sigma} \expp{- \langle
  \rho_s,f \rangle} \left[\psi( f(H_s)) - 
     f'(H_s) \ind_{\{H_s<\infty\}}\right]\; ds
\]
is a martingale w.r.t. the filtration generated
by $\rho$.  
\end{cor}

\begin{rem}
   \label{rem:H}
   Notice from Definition 1.2.1 in \cite{dlg:rtlpsbp} that $H_t<\infty $
   $\P_\mu$-a.s. for  all $t>0$ if  $\mu\in \cm_f$. Hence  the indicator
   $\ind_{\{H_s<\infty    \}}$   can    be   removed    in   Corollaries
   \ref{cor:mart1}  and  \ref{cor:mart-ts}  under $\P_\mu$  for  $\mu\in
   \cm_f$.
\end{rem}

\begin{rem}
  \label{rem:AD} {F}rom \cite{ad:falp}  (see Theorem 6.1 and Proposition
  6.3 in  Section 6.1),  there is a  relationship between  the infinitesimal
  generator of $\rho$ and the infinitesimal generator of the exploration
  process, $\rho^{(\theta)}$,  associated to a Lévy  process with exponent
  $\psi^{(\theta)}$, where
\[
\psi^{(\theta)}(\lambda)=\psi(\lambda+\theta) -\psi(\theta), \quad
\lambda\geq 0.
\]
More precisely, let  $F, K\in \cb(\cm_f)$  bounded  such that we
have    $\displaystyle    \E_\mu\left[\int_0^\sigma
  \val{K(\rho_s)}\;   ds  \right]<\infty   $   for any
$\mu\in     \cm_f$ and  $M_t=F(\rho_{t\wedge
  \sigma}) -\int_0^{t\wedge \sigma} K(\rho_s)$, for $t\geq 0$, defines a
martingale under $\P_\mu$. 
Then, we have for all $\mu\in \cm_f$, $\P_{\mu}$-a.s.
\begin{equation}
   \label{eq:majo-rq}
\int_0^{ \sigma^{(\theta)} } du \int_{(0,\infty )} \left(1-\expp{-\theta
    \ell}   \right)\;\pi(d\ell)\val{F([  \rho_u^{(\theta)},  \ell\delta_0])-
  F(\rho^{(\theta)}_u)} <\infty ,
\end{equation}
where $\sigma^{(\theta)}=\inf\{t>0;\rho^{(\theta)}_t=0\}$, 
and the process $(N_t, t\geq 0)$ is a martingale, where for $t\geq 0$, 
\[
N_t=F( \rho_{t\wedge  \sigma^{(\theta)}}^{(\theta)})  - 
\int_0^{t\wedge   \sigma^{(\theta)}} du\; \left( K( \rho_u^{(\theta)}) + 
\int_{(0,\infty )} \left(1-\expp{-\theta \ell}
\right)\;\pi(d\ell)\Big(F([ \rho_u^{(\theta)}, \ell\delta_0])- F(
\rho_u^{(\theta)})\Big)\right).
\]

Now, if we take $F$ and $K$ defined by
\reff{eq:def-FK}, and assume that $f\geq \varepsilon$ for some constant
$\varepsilon>0$. Then, we have 
\[
\E_\mu\left[\int_0^\sigma
  \val{K(\rho_s)}\;   ds  \right]
\leq  \E_\mu\left[\int_0^\sigma \expp{-\varepsilon \langle \rho_s, 1 \rangle}
  \;   ds  \right] 
= \frac{1 - \expp{- \varepsilon \langle \mu,1
    \rangle} }{\psi(\varepsilon)} <\infty .
\]
Notice that  $F([  \mu,  \ell\delta_0])-
  F(\mu)\leq 0$ and 
\begin{align*}
    \int_{(0,\infty )} \!\left(1-\expp{-\theta
    \ell}   \right)\;\pi(d\ell)\Big( F([  \mu,  \ell\delta_0])-
  F(\mu)\Big)  
&= - F(\mu) \int_{(0,\infty )}\! \left(1-\expp{-\theta
    \ell}   \right)\left(1-\expp{-
    \ell  f(H^\mu) }   \right) \;\pi(d\ell)\\
&=- F(\mu) \left(\psi(f(H^\mu)) - \psi^{(\theta)}(f(H^\mu))\right).
\end{align*}
Condition \reff{eq:majo-rq} is also satisfied. We get
for  the martingale $(N_t, t\geq 0)$ given above that 
\begin{align*}
N_t
&=F( \rho_{t\wedge  \sigma^{(\theta)}}^{(\theta)})  - 
\int_0^{t\wedge   \sigma^{(\theta)}} du\; \left( K( \rho_u^{(\theta)}) + 
\int_{(0,\infty )} \left(1-\expp{-\theta \ell}
\right)\;\pi(d\ell)\Big(F([ \rho_u^{(\theta)}, \ell\delta_0])- F(
\rho_u^{(\theta)})\Big)\right)    \\
&=F( \rho_{t\wedge \sigma^{(\theta)}}^{(\theta)})  - 
\int_0^{t\wedge   \sigma^{(\theta)}} du\; F( \rho_s) 
\left[\psi^{(\theta)} ( f(H^{\rho^{(\theta)}_s}) ) -
  f'(H^{\rho^{(\theta)}_s}) \right]. 
\end{align*}
And thanks to Remark \ref{rem:H}, we recover Corollary \ref{cor:mart-ts}
for the process $\rho^{(\theta)}$ for any initial measure in $\cm_f$.
\end{rem}

\section{Proof of Theorem \ref{th:feller}}
\label{sec:Feller}
In this sequel $x_+=\max(x,0)$ denotes the positive part of $x \in \R$. 

\subsection{A distance which induces the topology of weak convergence}
 
Let  $G$   be an  increasing  non-negative  bounded
function of  class $\cc^1$ on $\R_+$. For  convenience, we write $G(\infty  )$ for
the  limit of $G(t)$  as $t$  goes to  infinity, and  we shall  assume that
$G(0)=0$ and  $G(\infty )\leq 1$.  The set  $E=\R_+\cup\{\infty \}$
endowed with
the  distance $d(x,y)=|G(x)-G(y)|$ is  compact. Then  the set  of finite
measures  on  $E$, $\cm_f(E)$,  equipped  with  the  topology of weak convergence  is
locally compact.

For  $r\in \R_+$,  $\mu\in \cm_f(E)$, we  set $H_r^\mu=H^{k_r\mu}$.
Notice the function $r\mapsto H^\mu_r$ is non-increasing, right
continuous, satisfies $H^\mu_{\langle \mu,1 \rangle}=0$. For $r\in
[0,\langle  \mu,1 \rangle]$, we have   $H_r^\mu=g\big
((\langle \mu,1 \rangle  - r)^-\big)$, where $g$ is  the right continuous
inverse of the cumulative  distribution function of $\mu$. 
In particular, if for $\mu,\nu\in \cm_f(E)$ we have $\langle \mu,1
\rangle=\langle \nu,1 \rangle$ and $H^\mu_r=H^\nu_r$ for all $r\geq 0$,
then $\mu=\nu$. 

\begin{lem}
\label{lem:propH}
For all $r,v\in \R_+$, $\mu\in \cm_f(E)$, we have 
\begin{equation}
   \label{eq:v-Hr-equiv}
v<H_r^\mu \Longleftrightarrow \mu((v,+\infty ]) >r.
\end{equation}
For any $h\in \cb_b(E)$,  we have 
\begin{equation}
   \label{eq:hH=hmu1}
\int_0^{\langle \mu,1 \rangle}  h(H^\mu_r)\;dr=\langle \mu,h \rangle.
\end{equation}
and if  $h(0)=0$ this can be also written as 
\begin{equation}
   \label{eq:hH=hmu}
\int_0^\infty  h(H^\mu_r)\;dr=\langle \mu,h \rangle.
\end{equation}

\end{lem}

\begin{proof}
   Notice that $H^\mu_r=\sup\{ x\in E; \mu([x, +\infty ])>r\}$. This
   implies \reff{eq:v-Hr-equiv}. Let $h$ be a non-negative
   non-decreasing bounded
   function defined on $\R_+$ of class $\cc^1$ such that $h(0)=0$ and 
$\displaystyle   \int_{\R_+} \val{h'(v)} \mu([v,+\infty ]) dv $ is
finite. We have 
\begin{align*}
  \int_0^\infty  h(H^\mu_r)\;dr 
& = \int_{\R_+^2} h'(v) \ind_{\{v<H^\mu_r\}} \;dr dv\\
& = \int_{\R_+^2} h'(v) \ind_{\{r<\mu((v,+\infty ])\}} \;dr dv\\
&= \int_0^\infty  h'(v) \mu((v,+\infty ]) \; dv\\
&= \int_{\R_+\times E} h'(v) \ind_{\{u> v\}} \; dv\mu(du)\\
&= \langle \mu,h \rangle,
\end{align*}
where we used \reff{eq:v-Hr-equiv} for the second equality, and $h(0)=0$
for the first and last equality. By linearity and monotone class
Theorem, we get \reff{eq:hH=hmu} for any $h\in \cb_b(E)$, such that
$h(0)=0$. Then \reff{eq:hH=hmu1} is a direct consequence of
\reff{eq:hH=hmu} as $H^\mu_r=0$ for $r>\langle \mu,1 \rangle$. 
\end{proof}

In  particular,   \reff{eq:hH=hmu}  holds  for  $h=G$.    Thus  we  have
$\int_0^\infty G(H^\mu_r)dr \leq \langle  \mu,1 \rangle<\infty $ for all
$\mu\in  \cm_f(E)$.   We introduce  the  following  function defined  on
$\cm_f(E)^2$: for $\mu,\nu\in \cm_f(E)$
\begin{align*}
 D(\mu,\nu)
&=
 |\langle \mu,1 \rangle - \langle \nu,1 \rangle|+ \int_0^\infty
 d(H_r^\mu, H_r^\nu) \; dr \\
&= |\langle \mu,1 \rangle - \langle \nu,1 \rangle|+ \int_0^{\max(\langle
  \mu,1 \rangle ,\langle \nu,1 \rangle) }d(H_r^\mu,H_r^\nu)\; dr.
\end{align*}
Since $G$ is bounded by 1, we have
\begin{equation}
   \label{eq:majo-D}
\langle \mu,1 \rangle\leq  D(0,\mu)\leq  2 \langle \mu,1 \rangle. 
\end{equation}
{F}rom  what precedes  Lemma  \ref{lem:propH}, we  get that  $D(\mu,\nu)=0$
implies $\mu=\nu$.  It then is easy to check that $D$ is a distance on
$\cm_f(E)$. 

\begin{rem}
   \label{rem:suite-non-compacte}
   We deduce from \reff{eq:majo-D} that a sequence $(\mu_n, n\geq 1)$ is
   unbounded in  $(\cm_f(E), D)$ if  and only if the  sequence $(\langle
   \mu_n, 1 \rangle, n\geq 1)$ is unbounded.
\end{rem}

\begin{lem}
The topology induced by $D$ on $\cm_f(E)$ corresponds to the topology
induced by the weak convergence. 
\end{lem}

Notice that  $(\cm_f(E),D)$ is a locally compact Polish space.

\begin{proof} 

  We  first  consider a  sequence  $(\mu_n,  n\geq  1)$ of  elements  of
  $\cm_f(E)$ which  converges for the  distance $D$ to $\mu$.   We shall
  prove  that this  sequence  converges  weakly to  $\mu  $.  Let  $f\in
  \cb_b(E)$  such that  $\val{f(x)-f(y)}  \leq d(x,y)$  for all  $x,y\in
  E$. We set $f^*$ for $f -f(0)$.  Thanks to \reff{eq:hH=hmu}, we have
\begin{align*}
  |\langle\mu_n,f\rangle -\langle\mu,f\rangle | &=\Big| \langle\mu_n,
  f(0)\rangle  - \langle\mu,f(0)\rangle + \int_0^\infty  [f^*(H_r^{\mu_n})
  -f^*(H_r^\mu)] \; dr  \Big|\\ 
  &\leq  |f(0)||\langle\mu_n,1\rangle -\langle\mu,1\rangle | +\int_0^\infty 
  |f^*(H_r^{\mu_n}) - f^*(H_r^\mu)| \; dr  \\
  &\leq |f(0)||\langle\mu_n,1\rangle -\langle\mu,1\rangle | + \int_0^\infty 
  d(H_r^{\mu_n},H_r^\mu) \; dr  \\
  &\leq ( |f(0)|+1 ) D(\mu_n, \mu).
\end{align*}
As this holds for any Lipschitz function, we get that $(\mu_n, n\geq 1)
$ converges weakly to $\mu $.

Let $(\mu_n, n\geq 1)$ be a sequence of elements of $\cm_f(E)$ which
converges weakly to $\mu $.  In particular, the sequence
$(\langle\mu_n,1\rangle, n\geq 1  )$ converges to  $\langle\mu,1\rangle
$. Thus, there exists a finite constant $A$ such that
$\langle\mu,1\rangle\leq A$ and 
$\langle\mu_n,1\rangle\leq  A$ for all $n\geq 1$. Then, we have 
\begin{align*}
   \int_0^\infty d(H^{\mu_n}_r, H^\mu_r) \; dr
&=\int_0^\infty dr\; |G(H^{\mu_n}_r) -G(H^\mu_r)|\\
&=\int_0^\infty  dr\; \Big |\int_0^\infty dv \; G'(v) [\ind_{\{r<\mu_n((v,\infty ])\}} -
\ind_{\{r<\mu((v,\infty ])\}}]\Big|\\
&\leq  \int_{\R_+^2}   drdv \; G'(v) |\ind_{\{r<\mu_n((v,\infty ])\}} -
\ind_{\{r<\mu((v,\infty ])\}}|\\
&= \int_{\R_+}   dv \; G'(v) |\mu_n((v,\infty ])- \mu((v,\infty ])|,
\end{align*}
where we  used \reff{eq:v-Hr-equiv} for  the second equality.   The weak
convergence of $(\mu_n, n\geq  1)$ towards $\mu$ implies that $dv$-a.e.,
$\lim_{n\rightarrow  \infty }  \mu_n((v,\infty ])=\mu((v,\infty  ])$. As
$G'(v) |\mu_n((v,\infty ])- \mu((v,\infty ])|\leq 2 A G'(v)$ and $2AG' $
is integrable, we deduce  from dominated convergence Theorem that $\displaystyle
\lim_{n\rightarrow\infty }  \int_0^\infty \!\!\!d(H^{\mu_n}_r, H^\mu_r)\; dr=0
$. This and the convergence of $(\langle \mu_n, 1 \rangle, n\geq 1) $ to
$\langle  \mu,1   \rangle$  imply  that   $\displaystyle \lim_{n\rightarrow  \infty  }
D(\mu_n,\mu)=0$.
\end{proof}

We introduce the distance $D$ in order to get good continuity properties
for the  concatenation and  the erasing functions $k_a$.
\begin{lem}
Let  $a\geq 0$. We have   $D(k_a\mu, k_a\nu)\leq D(\mu,\nu)$.
\end{lem}

\begin{proof}
Notice $H^{k_a\rho}_r=H^\rho_{r+a}$.  We have
   \begin{align*}
      D(k_a\mu, k_a\nu)
&= |(\langle\mu,1\rangle -a)_+ -
(\langle\nu,1\rangle -a)_+ |+ \int_0^\infty  d(H^\mu_{r+a}, H^\nu_{r+a})\;dr \\
&= |(\langle\mu,1\rangle -a)_+ -
(\langle\nu,1\rangle -a)_+ |+ \int_a^\infty  d(H^\mu_{r}, H^\nu_{r})\;dr \\
&\leq D(\mu,\nu),
   \end{align*}
where we used for the last inequality the fact that  the function  $x\mapsto (x-a)_+$ is Lipschitz with constant  equal to one.
\end{proof}

With the convention
$x+\infty =\infty$  for all  $x\in E$, we have 
\[
H^{[k_a\mu,\rho]}_r=
\begin{cases}
     H^\mu_a+H^\rho_r   & \text{if $r\leq \langle\rho,1\rangle $,}\\
H^\mu_{a+r-\langle\rho,1\rangle } &\text{if $r\geq  \langle\rho,1\rangle $.}
\end{cases}
\]

Let $\mu\in \cm_f(E)$, and  $\ca_\mu\subset
\R_+$ be the set of continuity points of the function $r\mapsto
H^\mu_r$. Notice that $\R_+\backslash\ca_\mu$ is at most countable. 

\begin{lem}
\label{lem:erase-D}
Let $(\mu_n,  n\geq 1)$  be a sequence  of $\cm_f(E)$ which  converges
weakly to
$\mu\in \cm_f(E)$. For all $a\in \ca_\mu$,
$\rho\in \cm_f(E)$, we have that $([k_a\mu_{n},\rho], n\geq 1)$
converges weakly to $ [k_a\mu,\rho]$: 
\[
\lim_{n\rightarrow \infty } D([k_a\mu_{n},\rho], [k_a\mu,\rho])=
0.
\]
\end{lem}

\begin{proof}
  Let $(\mu_n,  n\geq 1)$ be a  sequence of elements  of $\cm_f(E)$ which
  converges  to  $\mu\in  \cm_f(E)$   for  the  distance  $D$  (i.e. which
  converges weakly).  
  For all $a\geq 0$,  we have $  \langle [k_a\mu_n,\rho], 1
  \rangle =\langle \rho,1 \rangle  + (\langle \mu_n,1 \rangle -a)_+$ and
  $\lim_{n\rightarrow  \infty  }   \langle  [k_a\mu_n,\rho],  1  \rangle
  =\langle [k_a\mu,\rho],  1 \rangle$.   There exists a  finite constant
  $A$      such      that      $\langle\mu,1\rangle\leq      A$      and
  $\langle\mu_n,1\rangle\leq  A$  for all  $n\geq  1$.  Furthermore  the
  sequence  of functions  $(G(H^{\mu_n}_\cdot), n\geq  1)$  converges in
  $L^1([0,A], dr)$  to $G(H^{\mu}_\cdot)$.  Since  $G(H^{\mu_n}_\cdot)$ and $G(H^{\mu}_\cdot)$ are
  non-increasing, we deduce that $\lim_{n\rightarrow \infty }
  G(H_a^{\mu_n}) =G(H_a^\mu) $ for any continuity  point $a$ of $G(H_ \cdot^\mu) $.
Since $G$ is  increasing and
  continuous, we deduce that 
for any $a\in \ca_\mu\cap[0,A]$, we
  have 
  $\lim_{n\rightarrow \infty } H_a^{\mu_n} =H_a^\mu$.
The results is also trivially true for $a>A$. 

Let $a\in \ca_\mu$. For $r\leq \langle\rho,1\rangle $, we have 
\[
H^{[k_a\mu_{n},\rho]}_r=H^{\mu_{n}}_a     +     H^\rho_r
\xrightarrow[n\rightarrow\infty
]{}
 H^{\mu}_a + H^\rho_r= H^{[k_a\mu,\rho]}_r,
\]
and for $r\geq \langle\rho,1\rangle $, $dr$-a.e.
\[
H^{[k_a\mu_{n},\rho]}_r=H^{\mu_{n}}_{a +r
  -\langle\rho,1\rangle } \xrightarrow[n\rightarrow\infty 
]{}
 H^{\mu}_{a
  +r -\langle\rho,1\rangle }= H^{[k_a\mu,\rho]}_r.
\]
Notice that for $r\geq  A+ \langle \rho,1 \rangle$, we have 
$ H^{[k_a\mu_{n},\rho]}_r=H^{[k_a\mu,\rho]}_r=0$. 
Since $G$ is continuous, we get by dominated convergence  for $a\in
\ca_\mu$, that 
\[
\lim_{n\rightarrow\infty } \int_0^\infty 
d(H^{[k_a\mu_{n},\rho]}_r, H^{[k_a\mu,\rho]}_r)\; dr =
0. 
\]
This    implies     that    for $a\in \ca_\mu$,  $\lim_{n\rightarrow\infty    }
D([k_a\mu_{n},\rho], [k_a\mu,\rho]) = 0$.
\end{proof}

Let $\cc_0(\cm_f(E))$ be the set of real continuous functions defined on
$(\cm_f(E),D)$  such  that  $\displaystyle \lim_{n\rightarrow  \infty  }
F(\mu_n)=0$   whenever    $\displaystyle   \lim_{n\rightarrow  \infty  }   D(0,
\mu_n)=+\infty  $, that  is  whenever $\displaystyle  \lim_{n\rightarrow
  \infty } \langle \mu_n, 1 \rangle=\infty $.

\begin{cor}\label{cor:continuity}
  Let   $(I^*,\rho)$    be   a   random   variable    with   values   in
  $\R_+,\times\cm_f(E)$, such  that the distribution  of $I^*$ has  no atom
  (i.e. $\P(I^*=x)=0$ for all $x\in \R_+$).  Let $F\in \cc_0(\cm_f(E))$.
  Then the function  $\cf:\mu\mapsto \E[F([k_{I^*} \mu, \rho])]$ belongs
  to $\cc_0(\cm_f(E))$.
\end{cor}

\begin{proof}
   Let $Q$ be the distribution of $I^*$
   and $(P_a(d\rho), a\geq 0)$ be a measurable  version of the
   conditional law of $\rho$ knowing $I^*$: $P_a(d\rho)\; Q(da) $ is the
   distribution of $(I^*,\rho)$. For $F\in \cc_0(\cm_f(E))$, we have 
\[
\E[F([k_{I^*} \mu, \rho])] =\int P_a[F([k_{a} \mu, \rho])]\; Q( da).
\]

Let $(\mu_n, n\leq 1)$ be  a sequence converging to $\mu$.  {F}rom Lemma
\ref{lem:erase-D},    for     all    $a\in    \ca_\mu$,    $\displaystyle
\lim_{n\rightarrow  \infty  }  D([k_a\mu_{n},\rho], [k_a\mu,\rho])=  0$.
Since  the complementary  of $\ca_\mu$  is at  most countable  and since
$Q(da)$  has   no  atoms,  we  get   that  $Q(da)$-a.s.   $\displaystyle
\lim_{n\rightarrow  \infty }  F([k_{a} \mu_{n},  \rho]) =  F([k_{a} \mu,
\rho])  $,  for any  $\rho\in  \cm_f(E)$.   Notice  $F$ is  bounded.  By
dominated  convergence, we  get  that $\displaystyle  \lim_{n\rightarrow
  \infty }  \E[F([k_{I^*} \mu_{n}, \rho])] =  \E[F([k_{I^*} \mu, \rho])]
$.  Thus the function $\cf$ is continuous.

Notice  that  the  total mass  of  $[k_{I^*}  \mu,  \rho]$ is  equal  to
$(\langle\mu,1\rangle -I^*)_++\langle\rho,1\rangle  $. If $\displaystyle
\lim_{n\rightarrow  \infty  }   \langle\mu_n,1\rangle  =\infty  $,  then
a.s. we have $\lim_{n\rightarrow \infty } F([k_{I^*} \mu_n,\rho])=0$. By the
dominated  convergence  Theorem,  we  get $\lim_{n\rightarrow  \infty  }
\cf(\mu_n)= 0$.

In conclusion, we get that  $\cf\in C_0(\cm_f(E))$. 
\end{proof}

\subsection{Feller property: Proof of Theorem \ref{th:feller}}

Let $(P_t, t\geq 0)$ be the transition semi-group of the Markov process
$\rho$ on $\cm_f(E)$. 

Recall $\cc_0(\cm_f(E))$ is the set of real continuous functions defined on
$(\cm_f(E),D)$  such  that  $\displaystyle \lim_{n\rightarrow  \infty  }
F(\mu_n)=0$   whenever    $\displaystyle   \lim_{n\rightarrow  \infty  }   D(0,
\mu_n)=+\infty  $, that  is  whenever $\displaystyle  \lim_{n\rightarrow
  \infty } \langle \mu_n, 1 \rangle=\infty $.

Let us recall
that (see e.g. \cite{ry:bm2}, Proposition III.2.4) $\rho$ is a
Feller process if and only if
\begin{itemize}
\item[(i)] If $F\in \cc _0(\cm_f(E))$, then  for
  every $t\ge 0$, $P_tF  \in \cc_0(\cm_f(E))$. 
\item[(ii)] For every $F\in\cc_0(\cm_f(E))$, for every $\mu\in\cm_f(E)$,
$\displaystyle\lim_{t\to 0}P_tF(\mu)=F(\mu)$.
\end{itemize}

To prove Condition (i), we state the next Lemma.

\begin{lem}\label{lem:atoms}
For every $t>0$, the distribution of $I_t$
has no atom.
\end{lem}

\begin{proof}
  Let us  denote by $\tau=(\tau_r,r\ge 0)$  the right-continuous inverse
  of the  process $(-I_t,t\ge  0)$. We now  (cf \cite{b:pl},  Chap. VII)
  that  the  process $\tau$  is  a  subordinator  with Laplace  exponent
  $\psi^{-1}$.        As        $\lim_{\lambda\rightarrow\infty        }
  \psi(\lambda)/\lambda=\infty      $,      we      have     that      $
  \lim_{\lambda\rightarrow\infty  }  \psi^{-1}(\lambda)/\lambda=0$.  The
  process $\tau$ has no drift.  It is moreover strictly increasing since
  the  process  $-I_t$ is  continuous.   It is  easy  to  check that  if
  $\P(-I_t=x)>0$ for some $x>0$, then we have $\P(-I_t=x, \tau_{x}=t)>0$
  and $\P(\tau_{-I_t}=t)>0$.  This is in contradiction to  the fact that
  $\tau$ has no  drift (see Theorem 4 in  \cite{b:pl}).  Hence the Lemma
  is proved.
\end{proof}

Let $F\in
\cc_0(\cm_f(E))$ and let $t\ge 0$.  Then,
as the distribution of the random variable $I_t$ has no atoms (see
Lemma \ref{lem:atoms}),
we can apply 
Corollary \ref{cor:continuity} with  $I^*=-I_t$ and get the function
$\cf: \mu \mapsto \E_\mu[F(\rho_t^\mu)]$
belongs to $\cc_0\left(\cm_f(E)\right)$.

It remains to  prove Condition (ii). Let us  remark that the exploration
process  is  right-continuous at  $t=0$.  Indeed,  for every  continuous
function $f$ on $E$ bounded by $1$, we have, $\P_\mu$-a.s.
\begin{align*}
\left|\langle \rho_t,f\rangle-\langle \mu,f\rangle\right| &
=\left|\left\langle\left[k_{-I_t}\mu,\rho^0_t\right],f\right\rangle-\langle
\mu,f\rangle\right|\\
& =\left|- \left\langle
 \mu - k_{-I_t}\mu,f\right\rangle+\left\langle\rho_t^0,f(H_{-I_t}^\mu+\cdot)\right\rangle\right|\\
& \le \left\langle
\mu-k_{-I_t}\mu,1\right\rangle+\langle\rho_t^0,1\rangle\\
&\leq -I_t+X_t-I_t.
\end{align*}
And the right continuity of $\rho$ at 0 follows from the
right-continuity of $X$ and $I$ at 0. 
Let $F\in\cc_0(\cm_f(E))$. As the function $F$ is bounded, we
get by dominated convergence that, for every $\mu\in\cm_f(E)$,
$$\E_\mu[F(\rho_t)]\xrightarrow[t\to 0]{} F(\mu).$$

\section{Proof of Theorem \ref{th:inf-gene-0}}
\label{sec:inf-gene}

\subsection{Other properties of the exploration process}
In this section we recall some properties of the exploration process. 
\subsubsection{Poisson decomposition}
 \label{sec:rho-Poisson}
 Let $\mu\in  \cm_f(E)$. We decompose  the path of $\rho$  under $\P_{\mu}$
 according to excursions of the  total mass of $\rho$ above its minimum.
 More precisely, we  denote by  $(\alpha_i,\beta_i),  i\in \ci$  the
 excursion  intervals of  $X-I$ away  from $0$.   For every  $i\in \ci$,
 $t\in          (\alpha_i,          \beta_i)$,          we          have
 $\rho_{\alpha_i}=k_{-I_{\alpha_i}}\mu=k_{-I_{\alpha_i}}   \rho_t$.   We
 define $\rho^i$ by the formula $\rho^i_t=\rho^0_{(\alpha_i +t )\wedge
   \beta_i}$ or equivalently if $\mu\in \cm_f$:
\[
[k_{-I_{\alpha_i}}\mu, \rho^i_t]=\rho_{(\alpha_i +t )\wedge \beta_i}.
\]
The local time of $X-I$ at level  $0$ is given by $-I$, and if $(\tau_r,
r\geq  0)$ is  the right  continuous inverse  of $I$,  then  the process
$(\tau_r,r\ge 0)$  is a  subordinator with Laplace  exponent $\psi^{-1}$
(cf  \cite{b:pl}, Chap. VII).   Excursion theory  for $X-I$  ensures the
following result.
\begin{lem}
\label{lem:rho-decomp}
   Let $\mu\in \cm_f(E)$. The point
   measure $\displaystyle \sum_{i\in I} \delta_{(-I_ {\alpha_i},\rho^i)}$ is under
    $\P_{\mu}$ a Poisson point measure with intensity $dr \N[d\rho]$.

Let $F $ be a non-negative measurable function
   defined on $\R_+\times \cm_f(E)\times \mathbb{D}([0,\infty ),\cm_f(E))$, where
   $\mathbb{D}([0,\infty ),\cm_f(E))$ stands for the Skorohod space of càd-làg
   path in $\cm_f(E)$. We have 
\begin{equation}
   \label{eq:rho-decomp}
\E_\mu\left[\sum_{i\in \ci} F(\alpha_i,\rho_{\alpha_i}, \rho^i) \right] 
=\E_\mu \left[\int_0^\infty  dr\; \N\left[F(v,k_r\mu,\rho)
  \right]_{|v=\tau_r}\right]. 
\end{equation}
\end{lem}

\subsubsection{The dual process and representation formula}
\label{sec:dual}

We  shall need the  $\cm_f$-valued process  $\eta=(\eta_t,t\ge 0)$
 defined under $\N$ by
\[
\eta_t(dr)=\sum_{\stackrel{0<s\le t}{X_{s-}<I_t^s}}(X_s-I_t^s)\delta
_{H_s}(dr).
\]
The  process  $\eta$ is  the  dual process  of  $\rho$  under $\N$  (see
Corollary   3.1.6   in   \cite{dlg:rtlpsbp}).   Let   $\sigma=\inf\{s>0;
\rho_s=0\}$ denote the length of  the excursion.  Using this duality, it
easy to check (see the proof of Lemma 3.2.2 in \cite{dlg:rtlpsbp}), that
for any bounded measurable function $F$ defined on $\cm_f$,
\begin{equation}
   \label{eq:N-t-rho}
\N\left[\int_0^\sigma \expp{- \psi(\gamma) t } F(\rho_t)\right]
= \N\left[\int_0^\sigma \expp{- \gamma\langle \eta_t,1  \rangle  }
  F(\rho_t)\right] . 
\end{equation}

We recall the Poisson representation of $(\rho,\eta)$ under $\N$. Let
$\mathcal{N}(dx\,   d\ell\,  du)$   be  a   Poisson  point   measure  on
$[0,+\infty)^3$ with intensity
$$dx\,\ell\pi(d\ell)\ind_{[0,1]}(u)du.$$
For every $a>0$, let us denote by $\mathbb{M}_a$ the law of the pair
$(\mu_a,\nu_a)$ of finite measures on $\R_+$ defined by:  for $f\in \cb_+(\R_+)$ 
\begin{align}
\label{def:mu_a}
\langle \mu_a,f\rangle  & =\int\mathcal{N}(dx\, d\ell\, du)\ind_{[0,a]}(x)u\ell f(x),\\
\label{def:nu_a}
\langle \nu_a,f\rangle  & =\int\mathcal{N}(dx\, d\ell\, du)\ind_{[0,a]}(x)\ell(1-u)f(x).
\end{align}
We eventually set $\mathbb{M}=\int_0^{+\infty}da\, \expp{-\alpha_0
  a}\mathbb{M}_a$. 

\begin{prop}\label{prop:poisson_representation1}
For every non-negative measurable function $F$ on $\cm_f^2$,
\[
\N\left[\int_0^\sigma F(\rho_t, \eta_t) \; dt \right]=\int\mathbb{M}(d\mu\,
    d\nu)F (\mu, \nu) .
\]
\end{prop}

\subsection{Computation of the resolvent}

Let  $f$ be a  bounded function  defined on  $\R_+$ of  class  $\cc^1$ with
bounded  first derivative  such  that $f(\infty  )=\lim_{t\rightarrow\infty  }
f(t)$  exists. By convention, we put $f'(\infty )=0$. 

We set
for $\mu\in \cm_f(E)$
\begin{equation}
   \label{eq:def-Fl}
F_\lambda(\mu)=\expp{- \langle \mu,f \rangle} \left(\lambda -
  \psi(f(H^\mu)) + f'(H^\mu)\right).
\end{equation}
Let  $\gamma=\psi^{-1}(\lambda)$. We shall   compute 
\[
U_\lambda(F_\lambda)(\mu)=\E_\mu \left[\int_0^\infty dt\;
  \expp{- \psi(\gamma)  t - \langle \rho_t, f
    \rangle} \left(\psi(\gamma) -\psi(f(H_t)) + f'(H_t) \right) \right].
\]
We define  the function $f_r$ by  $f_r(\cdot)=f(\cdot+ H^\mu_r)$, where
$H^\mu_r=H^{k_r\mu}$ and use the convention $x+\infty =x$ for all $x\in
E$.  Using  notations of Section \ref{sec:rho-Poisson}
and \reff{eq:rho-decomp}, we have
\begin{align*}
  U_\lambda(F_\lambda)(\mu) 
&=\E_\mu\left[\!\sum_{i\in \ci} \expp{-\psi(\gamma) \alpha_i }
\int_0^{\beta_i-\alpha_i} \!\!\!\!dt\; \expp{-\psi(\gamma) t - \langle
  [k_{-I_{\alpha_i}}  \mu, \rho^i_t], f  \rangle}\right.\\
&\qquad\qquad\qquad\qquad\qquad\left.
\left(\psi(\gamma) -\psi(f_{-I_{\alpha_i}}(H^{\rho^i_t})) +
  f'_{-I_{\alpha_i}} (H^{\rho^i_t}) \right)\right]\\
  &= \E_\mu \left[\int_0^{\infty } \!\!dr\;  \expp{-
    \psi(\gamma)\tau_r} \N\left[\int_0^\sigma \!dt\; \expp{-\psi(\gamma) t
    - \langle [k_r \mu, \rho_t], f  \rangle}
\left(\psi(\gamma) -\psi(f_r(H_t)) + f'_r(H_t) \right) \right]
\right]. 
\end{align*}
Recall  $(\tau_r, r\geq 0)$ is a subordinator with Laplace exponent
$\psi^{-1}$ and \reff{eq:N-t-rho}.  Since $\langle [k_r
      \mu,\rho_t] , f \rangle= \langle k_r
      \mu , f \rangle + \langle \rho_t, f_r \rangle$, we get 
\[
   U_\lambda(F_\lambda)(\mu) 
=\int_0^{\infty } dr\;  \expp{-
    \gamma r} \expp{-\langle k_r \mu, f \rangle} \N\left[\int_0^\sigma
    dt\; \expp{-\gamma \langle \eta_t,1 \rangle - \langle \rho_t, f_r  \rangle}
\left(\psi(\gamma) -\psi(f_r(H_t)) + f'_r(H_t) \right) \right].
\]
We deduce from Proposition \ref{prop:poisson_representation1} that
\begin{multline*}
   \N\left[\int_0^\sigma
    dt\; \expp{-\gamma \langle \eta_t,1 \rangle - \langle \rho_t, f_r  \rangle}
\left(\psi(\gamma) -\psi(f_r(H_t)) + f'_r(H_t) \right) \right]\\
\begin{aligned}
 &  =\int_0^\infty  da\; \expp{- \alpha_0 a} \left(\psi(\gamma) - \psi(f_r(a)) +
  f'_r(a) \right)\\
&\hspace{2cm}
\exp\left\{- \int_0^a dx \int_0^1 du \int_{(0,\infty )} \ell\pi(d\ell)
    \left[1- \expp{- \gamma (1-u) \ell - u f_r(x)\ell} \right]\right\}\\
 &  =\int_0^\infty  da\; \left(\psi(\gamma) - \psi(f(a+H^\mu_r)) +
  f'(a+H^\mu_r) \right)
\expp{- \int_0^a dx \int_0^1 du\; \psi'(\gamma (1-u) + u
  f(x+H^\mu_r)) }.
\end{aligned}
\end{multline*}
We set for $y\in E$
\[
\Lambda(y)=\int_0^\infty  da\; \left(\psi(\gamma)  - \psi(f(a+y)) +
  f'(a+y) \right)
\expp{- g(a,y)},
\]
where $\displaystyle g(a,y)= \int_0^a dx\; \frac{\psi(f(x+y) ) -
  \psi(\gamma)}{f(x+y) - \gamma}   $ with the convention $\displaystyle \frac{\psi(v) -\psi(v)}{v-v}=\psi'(v)$,  so that 
\begin{equation}
   \label{eq:UlFlG}
 U_\lambda(F_\lambda)(\mu) =\int_0^{\infty } dr\;  \expp{-
    \gamma r} \expp{-\langle k_r \mu, f \rangle}  \Lambda(H^\mu_r).
\end{equation}
Since $\psi $ is positive increasing continuous and $f$ bounded, we have
$\lim_{a\rightarrow \infty } g(a,y)=\infty $.
Notice that for $y\in [0,\infty )$
\begin{align*}
   \Lambda(y)
&=\int_0^\infty  da\; \left(- \partial_a g(a,y)(f(a+y) - \gamma)  + f'(a+y)\right)
\expp{- g(a,y)}\\
&=\left[  ( f(a+y) - \gamma)  \expp{- g(a,y)}      \right]_0^\infty \\
&= \gamma-f(y).
\end{align*}
We also have $\Lambda(\infty )=\gamma-f(\infty )$.
We deduce from \reff{eq:UlFlG} that 
\[
U_\lambda(F_\lambda)(\mu) =\int_0^{\infty } dr\;  \expp{-
    \gamma r-\langle k_r \mu, f \rangle}  (\gamma - f(H^\mu_r)).
\]
Notice that $\langle k_r\mu,f \rangle=0 $ and $H^\mu_r=0$ for $r>
\langle \mu,1 \rangle $, so that we have 
\[
\int_{\langle \mu,1 \rangle}^{\infty } dr\;  \expp{-
    \gamma r-\langle k_r \mu, f \rangle}  (\gamma - f(H^\mu_r))
= \int_{\langle \mu,1 \rangle}^{\infty } dr\;  \expp{-
    \gamma r}  (\gamma - f(0))= \expp{- \gamma \langle \mu,1 \rangle}
  \left(1 -\frac{f(0)}{\gamma} \right).
\]
We deduce from \reff{eq:hH=hmu1} with $\mu$ replaced by $k_r\mu$ that
for $r\leq \langle \mu,1 \rangle$, 
\[
\langle k_r \mu, f \rangle=\int_0^{\langle \mu,1 \rangle -r}
f(H^{k_r \mu}_v)\; dv 
=\int_r^{\langle \mu,1 \rangle}  f(H^\mu_v)\; dv. 
\]
This implies that 
\begin{align*}
\int_0^{\langle \mu,1 \rangle} dr\;  \expp{-
    \gamma r-\langle k_r \mu, f \rangle}  (\gamma - f(H^\mu_r))
&=   \int_0^{\langle \mu,1 \rangle} dr\;  \expp{-
    \gamma r-\int_r^{\langle \mu,1 \rangle}  f(H^\mu_v)\; dv}  (\gamma -
  f(H^\mu_r))\\ 
&= \left[-\expp{-
    \gamma r-\int_r^{\langle \mu,1 \rangle}  f(H^\mu_v)\; dv}
\right]_0^{\langle \mu,1 \rangle}\\ 
&= \expp{- \langle \mu,f \rangle}- \expp{-\gamma {\langle \mu,1 \rangle}},
\end{align*}
where we used \reff{eq:hH=hmu1} for the first term of the right hand side
of the last equation. Eventually, we have
\[
U_\lambda(F_\lambda)(\mu) =\expp{- \langle \mu,f \rangle} -
\frac{f(0)}{\gamma} \expp{- \gamma \langle \mu,1 \rangle}. 
\]

\bibliographystyle{abbrv}
\bibliography{/home/delmas/cermics/Bibliographie/delmas}

\newcommand{\sortnoop}[1]{}
\begin{thebibliography}{1}

\bibitem{ad:falp}
R.~ABRAHAM and J.-F. DELMAS.
\newblock Fragmentation associated to {L}évy processes.
\newblock {\em Preprint CERMICS}, 2005.

\bibitem{a:crt3}
D.~ALDOUS.
\newblock The continuum random tree {III}.
\newblock {\em Ann. Probab.}, 21(1):248--289, 1993.

\bibitem{b:pl}
J.~BERTOIN.
\newblock {\em L\'evy processes}.
\newblock Cambridge University Press, Cambridge, 1996.

\bibitem{dlg:rtlpsbp}
T.~DUQUESNE and J.-F. LE~GALL.
\newblock {\em Random trees, {L}évy processes and spatial branching processes},
  volume 281.
\newblock Astérisque, 2002.

\bibitem{ek:mp}
S.~N. ETHIER and T.~G. KURTZ.
\newblock {\em Markov processes}.
\newblock Wiley, 1986.

\bibitem{lglj:bplplfss}
J.-F. LE~GALL and Y.~LE~JAN.
\newblock Branching processes in {L}évy processes: Laplace functionals of snake
  and superprocesses.
\newblock {\em Ann. Probab.}, 26:1407--1432, 1998.

\bibitem{lglj:bplpep}
J.-F. LE~GALL and Y.~LE~JAN.
\newblock Branching processes in {L}évy processes: The exploration process.
\newblock {\em Ann. Probab.}, 26:213--252, 1998.

\bibitem{ry:bm2}
D.~REVUZ and M.~YOR.
\newblock {\em Continuous martingales and {B}rownian motion}.
\newblock Springer Verlag, 2nd edition, 1995.

\end{thebibliography}

\end{document}